\documentclass{amsart}

\newtheorem{theorem}{Theorem}[section]
\newtheorem{lemma}[theorem]{Lemma}
\newtheorem{corollary}[theorem]{Corollary}

\theoremstyle{definition}
\newtheorem{definition}[theorem]{Definition}

\begin{document}

\title[Higher connectivity of ${\tt Hom}(G,K_n)$]
{A short proof of a conjecture
on the higher connectivity of graph
coloring complexes}

\author{Alexander Engstr\"om}

\address{ETH Z\"urich, IFW B29, 8092 Z\"urich, 
Switzerland} 
\email{engstroa@inf.ethz.ch}

\thanks{Research supported by ETH and Swiss 
National Science Foundation Grant PP002-102738/1}

\subjclass{57M15, 05C15}

\date{May 22, 2005}

\keywords{Graph homomorphisms, $k$-connectivity, 
{\tt Hom}--complexes}

\begin{abstract}
The ${\tt Hom}$--complexes were introduced
by Lov\'asz to study topological 
obstructions to graph colorings.
It was conjectured by Babson and Kozlov,
and proved by \v{C}uki\'c and Kozlov, that 
${\tt Hom}(G,K_n)$ is $(n-d-2)$--connected, 
where $d$ is the maximal degree of a vertex 
of $G$. We give a short proof of the
conjecture. 
\end{abstract}

\maketitle

\section*{Introduction}
It was conjectured by Babson and Kozlov 
\cite{BK}, and proved by \v{C}uki\'c and Kozlov 
\cite{CK}, that ${\tt Hom}(G,K_n)$ is 
$(n-d-2)$--connected, where $d$ is the maximal 
degree of a  vertex of $G$. We give a shorter 
proof of this, by generalizing the proof of that 
${\tt Hom}(K_m,K_n)$ is $(n-m-1)$--connected in 
Babson and Kozlov \cite{BK}.

For defintions and basic theorems on 
{\tt Hom}--complexes used in this text, see the 
papers mentioned above, or the survey by Kozlov 
\cite{K}.

\section{An analogue of the chromatic number}
An independent subset of vertices of a graph is 
a set, such that no vertices of it are adjacent. 
The minimal number of sets needed to partition 
the vertex set of a graph $G$ into indepentent 
sets is the chromatic number $\chi(G)$.

\begin{definition}
A \emph{covering} $I_1,I_2,\ldots, I_k$ of $G$ is
a sequence of independent subsets of $V(G)$ such 
that they partition $V(G)$, and $I_i$ is a 
maximal independent set in the induced subgraph 
of $G$ with vertex set $I_i \cup I_{i+1}\cup 
\ldots \cup I_k$, for all $i$, where $1\leq i 
\leq k$.
\end{definition}

A partition of $G$ into $\chi(G)$ independent 
sets can always be transformed to a covering by 
ordering the independent sets and if needed 
enlarging them. But a covering can use more than 
$\chi(G)$ sets. Define $\dot{\chi}(G)$ to be
the maximal number of sets in a covering of $G$.
Clearly $\dot{\chi}(G)\geq\chi(G)$.

\begin{lemma} \label{lemma:maxDegree}
If $d$ is the maximal degree of a vertex of $G$,
then $\dot{\chi}(G)\leq d+1$.
\end{lemma}
\begin{proof}
Let $I_1,I_2,\ldots,I_{\dot{\chi}(G)}$ be a 
covering of $G$, and $v\in I_{\dot{\chi}(G)}$. 
For each $i$, where $1\leq i < \dot{\chi}(G)$, 
there is a $w\in I_i$ adjacent to $v$, because 
otherwise $I_i$ would not be a maximal 
independent set. Hence the degree of $v$ is at 
least $\dot{\chi}(G)-1$. The degree of $v$ is at 
most $d$, thus $\dot{\chi}(G)\leq d+1$.
\end{proof}

\begin{lemma} \label{lemma:removeSome}
If $H$ is an induced subgraph of $G$, then
$\dot{\chi}(H)\leq \dot{\chi}(G)$.
\end{lemma}
\begin{proof}
It suffices to prove this when $H$ and $G$ only 
differ by a vertex $v$ of $G$. Let $I_1,I_2,
\ldots,I_{\dot{\chi}(H)}$ be a covering of $H$. 
If $v$ is adjacent to a vertex in each of the 
sets $I_i$, then $\{v\},I_1,I_2,\ldots, 
I_{\dot{\chi}(H)}$ is a covering of $G$ and 
$\dot{\chi}(H)+1 \leq \dot{\chi}(G)$. Otherwise,
let $I_j$ be the first set in the covering such 
that $v$ is not adjacent to any vertex of $I_j$. 
Then $I_1,I_2,\ldots, I_{j}\cup\{v\}, \ldots 
I_{\dot{\chi}(H)}$ is a covering of $G$, and
$\dot{\chi}(H)\leq \dot{\chi}(G)$.
\end{proof}

\begin{lemma} \label{lemma:removeMaxInd}
If $I$ is a maximal independent set of $G$, then
$\dot{\chi}(G)>\dot{\chi}(G\setminus I)$.
\end{lemma}
\begin{proof}
Let $I_1, I_2, \ldots , I_{\dot{\chi}(G\setminus
I)}$ be a covering of $G\setminus I$. Then $I,
I_1, I_2, \ldots , I_{\dot{\chi}(G\setminus I)}$
is a covering of $G$ with $1+\dot{\chi}
(G\setminus I)$ sets.
\end{proof}

\section{Higher connectivity of 
${\tt Hom}(G,K_n)$}

\begin{lemma}\label{lemma:collapse}
If $I$ is an independent set of $G$, and $I'
\subset I$, then $\Delta=\{\eta\in {\tt Hom}
(G,K_n) | n\in \eta(i)\Rightarrow i\in I\}$ 
collapses onto $\Delta'= \{\eta\in {\tt Hom}
(G\setminus (I \setminus I'),K_n) | n\in \eta(i)
\Rightarrow i\in I'\}$.
\end{lemma}
\begin{proof}
It suffices to prove this when $I\setminus I'
=\{v\}$. Let $\eta_1, \eta_2, \ldots \eta_k$ be 
an ordering of $\{\eta\in \Delta | n \not\in 
\eta(v)\}$ such that if $\eta(w) \supseteq 
\eta'(w)$ for all $w\in V(G)$ then $\eta$ is not
after $\eta'$. Define $\eta_i^\ast$ as 
$\eta_i^\ast(w)=\eta_i(w)$ for $w\neq v$, and 
$\eta_i^\ast(v)=\eta_i(v)\cup\{n\}$. Each 
successive removal of $\eta_i^\ast$ together with
$\eta_i$ from $\Delta$ for $i=1,2,\dots k$ is a 
collapse step. The cells left are $\Delta''=
\{\eta \in \Delta | \eta(v)=\{n\}\}$. Finally, 
there is a bijection between the face posets of 
$\Delta'$ and $\Delta''$ by extending each $\eta
\in \Delta'$ with $\eta(v)=\{n\}$.
\end{proof}

The main use of lemma \ref{lemma:collapse} is 
when $I'=\emptyset$. Then $n\not\in \eta(w)$ for
all $\eta\in \Delta'$ and $w\in V(G)\setminus I$,
so $\Delta'= {\tt Hom}(G\setminus I,K_{n-1})$.
Another way to prove the lemma is to use discrete
Morse theory \cite{F}.

We will use a variation of a Nerve Lemma, 
(Bj\"orner 10.6(ii) \cite{B}, Bj\"orner et.al. 
\cite{BLVZ}). A regular cell complex $\Delta$ is 
$m$--connected if there is a family of 
subcomplexes $\{ \Delta_i \}$ such that $\Delta=
\cup \Delta_i$, all of the subcomplexes 
$\Delta_i$ are $m$--connected, and all of the
intersection of several $\Delta_i$'s are 
$(m-1)$--connected.

\begin{theorem}
${\tt Hom}(G,K_n)$ is 
$(n-\dot{\chi}(G)-1)$--connected.
\end{theorem}
\begin{proof}
We use induction on $\dot{\chi}(G)$ and on 
$n-\dot{\chi}(G)$. When $\dot{\chi}(G)=1$, $G$ 
have no edges, so ${\tt Hom}(G,K_n)$ is 
contractible, and in particular 
$(n-\dot{\chi}(G)-1)$--connected. If 
$n-\dot{\chi}(G)=0$ then $n\geq \chi(G)$ so 
${\tt Hom}(G,K_n)$ is non-empty, and
$(n-\dot{\chi}(G)-1)$--connected.

For all $I\in \mathcal{I}$, let $\Delta_I=
\{\eta\in {\tt Hom} (G,K_n) | n\in \eta(i)
\Rightarrow i\in I\}$, where $\mathcal{I}$ is the
family of maximal independent subsets of $G$. 
Clearly ${\tt Hom}(G,K_n)=\cup_{I\in \mathcal{I}}
\Delta_I$. By lemma \ref{lemma:collapse}, the
complex $\Delta_I$ is homotopy equivalent to
${\tt Hom}(G\setminus I,K_{n-1})$, which is
$((n-1)-(\dot{\chi}(G)-1)-1)$--connected by
lemma \ref{lemma:removeMaxInd} and induction.
If $\mathcal{I}\supseteq \mathcal{I'}\neq 
\emptyset$ then $\cap_{I\in\mathcal{I'}} 
\Delta_I= \{\eta\in {\tt Hom}
(G,K_n) | n\in \eta(i)\Rightarrow i\in 
\cap_{I\in\mathcal{I'}} I \}$ is 
homotopy equivalent to ${\tt Hom} (G\setminus ( 
\cap_{I\in\mathcal{I'}} I ),K_{n-1})$
by lemma \ref{lemma:collapse}, and 
$((n-1)-\dot{\chi}(G)-1)$--connected by lemma 
\ref{lemma:removeSome} and induction. By the 
Nerve Lemma we are done.
\end{proof}

\begin{corollary}
${\tt Hom}(G,K_n)$ is $(n-d-2)$--connected.
\end{corollary}
\begin{proof}
Lemma \ref{lemma:maxDegree} states that 
$\dot{\chi}(G)\leq d+1$.    
\end{proof}

\bibliographystyle{amsplain}

\end{document}